\documentclass[12pt, a4paper, twoside]{article}
\usepackage{amssymb}
\usepackage[leqno]{amsmath}
\usepackage{latexsym}
\usepackage{enumerate}
\usepackage{amsmath, amssymb, amscd}
\Roman{equation}

\DeclareMathOperator{\Pic}{Pic}

\begin{document}

\title{\textbf{\Large{The characterization of aCM line bundles on quintic hypersurfaces in $\mathbb{P}^3$}}}

\author{Kenta Watanabe \thanks{Nihon University, College of Science and Technology,   7-24-1 Narashinodai Funabashi city Chiba 274-8501 Japan , {\it E-mail address:watanabe.kenta@nihon-u.ac.jp}, Telephone numbers: 090-9777-1974} }

\date{}

\maketitle 


\begin{abstract}

\noindent Let $X$ be a smooth quintic hypersurface in $\mathbb{P}^3$, let $C$ be a smooth hyperplane section of $X$, and let $H=\mathcal{O}_X(C)$. In this paper, we give a necessary and sufficient condition for the line bundle given by a non-zero effective divisor on $X$ to be initialized and aCM with respect to $H$. 

\end{abstract}

\section{Introduction} Let $X$ be a projective manifold of dimension $n$ over the complex number field $\mathbb{C}$, and let $H$ be a very ample line bundle on $X$. Let $E$ be a vector bundle on $X$. Then we call $E$ an {\it{arithmetically Cohen-Macaulay}} ({{\it{aCM}} for short) bundle with respect to $H$ if $H^i(E\otimes H^{\otimes l})=0$ for all integers $l\in\mathbb{Z}$ and $1\leq i\leq n-1$. Moreover, we say that $E$ is {\it{initialized}} if $E$ satisfies $h^0(E)\neq0$ and $h^0(E\otimes H^{\vee})=0$. If $X=\mathbb{P}^{n}$, an aCM bundle is characterized as a vector bundle obtained by a direct sum of line bundles on $X$ ([7]). However, this criterion is not correct for more general polarized manifolds. A vector bundle $E$ on $X$ is aCM with respect to $H$ if and only if $E\otimes H$ is so. Hence, previously, many people have studied indecomposable initialized aCM bundles on $X$ with respect to a given polarization $H$ and families of them. In particular, several results concerning initialized aCM bundles on smooth hypersurfaces in $\mathbb{P}^{n+1}$ are known. For example, if $X$ is a smooth quadric hypersurface in $\mathbb{P}^{n+1}$, then any non-split aCM bundle on $X$ of rank $r\geq2$ is isomorphic to a spinor bundle up to twisting by the hyperplane class of $X$ ([8]). In particular, if $n=2$, then any aCM bundle on $X$ splits. If $X$ is a smooth cubic hypersurface in $\mathbb{P}^3$,  Casanellas and Hartshorne have studied the families of initialized aCM bundles on $X$ ([2]). Moreover, Faenzi has given a complete classification of indecomposable  initialized aCM bundles of rank 2 on $X$ ([6]). If $X$ is a smooth quartic hypersurface in $\mathbb{P}^3$, Coskun and Kulkarni have constructed a 14-dimensional family of simple Ulrich bundles on $X$ of rank 2 with $c_1=H^{\otimes3}$ and $c_2=14$, in the case where $X$ is a Pfaffian quartic surface ([5, Theorem 1.1]). Casnati has classified indecomposable initialized aCM bundles of rank 2 on $X$, in the case where $X$ is general determinantal ([3]). 

In general, if an aCM bundle on $X$ splits into a direct sum of line bundles on $X$, then the line bundles on $X$ which appear in the splitting are also aCM. Conversely, a vector bundle on $X$ which is given by an extension of aCM line bundles on $X$ is aCM. Hence, the classification of aCM line bundles on $X$ is useful for constructing indecomposable aCM bundles on $X$ of higher rank. Indeed, Pons-Llopis and Tonini have classified aCM line bundles on a DelPezzo surface $X$ with respect to the anti-canonical line bundle on $X$, and constructed families of strictly semi-stable aCM bundles on $X$ ([9]). On the other hand, our previous work about the characterization of aCM line bundles on smooth quartic hypersurfaces in $\mathbb{P}^3$ ([11]) is deeply connected with the Casnati's work ([3]). In addition, recently, several other results (for example [4] and [12]) concerning the classification of aCM line bundles on polarized surfaces are also known. Since a smooth quartic hypersurface in $\mathbb{P}^3$ is a K3 surface and the properties of linear systems on K3 surfaces are well known ([10]), we can obtain the results as in [11] and [12] with comparative ease. On the other hand, if $X$ is a Del Pezzo surface or a ruled surface, any aCM line bundle on $X$ can be precisely denoted by using the generators of the Picard group of  $X$. However, in general, it is difficult to investigate initialized aCM line bundles on polarized surfaces. In this paper, we are interested in the characterization of aCM line bundles on smooth quintic hypersurfaces in $\mathbb{P}^3$. Our main theorem is as follows.

\newtheorem{thm}{Theorem}[section]

\begin{thm} 

Let $X$ be a smooth quintic hypersurface in $\mathbb{P}^3$, let $H$ be the hyperplane class of $X$, and let $C$ be a smooth member of the linear system $|H|$. Let $D$ be a non-zero effective divisor on $X$ of arithmetic genus $P_a(D)$, and let $k=C.D+1-P_a(D)$. Then $\mathcal{O}_X(D)$ is aCM and initialized if and only if the following conditions are satisfied.

$\;$

\noindent {\rm{(i)}} $0\leq k\leq 4$.

\smallskip

\smallskip

\noindent {\rm{(ii)}} If $0\leq k\leq 1$, then $C.D=10-k$ and $h^0(\mathcal{O}_C(D-C))=0$. 

\smallskip

\smallskip

\noindent {\rm{(iii)}} If $k=2$, then the following conditions are satisfied.

\smallskip

\smallskip

{\rm{(a)}} $C.D=1$ or $4\leq C.D\leq 8$.

{\rm{(b)}} If $C.D=7$, then $h^0(\mathcal{O}_X(2C-D))=1$.

{\rm{(c)}} If $C.D=8$, then $h^0(\mathcal{O}_C(D-C))=0$ and $h^0(\mathcal{O}_C(D))=3$.

\smallskip

\smallskip

\noindent {\rm{(iv)}} If $3\leq k\leq 4$, then the following conditions are satisfied.

\smallskip

{\rm{(a)}} $k-1\leq C.D\leq 10-k$.

{\rm{(b)}} If $8-k\leq C.D\leq 10-k$, then $h^0(\mathcal{O}_C(D))=5-k$. \end{thm} 

\noindent Our plan of this paper is as follows. In section 2, we give several basic results about line bundles on a smooth quintic hypersurface in $\mathbb{P}^3$. In section 3, we recall a known fact about aCM bundles on smooth quintic hypersurfaces  in $\mathbb{P}^3$ and give a lemma to prove Theorem 1.1. In section 4, we prove Theorem 1.1. 

$\;$

{\bf{Notations and conventions}}. In this paper, a curve and a surface are smooth and projective. Let $X$ be a curve or a surface. Then we denote the canonical bundle of $X$ by $K_X$. For a divisor or a line bundle $L$ on $X$, we denote by $|L|$ the linear system of $L$, and denote the dual of a line bundle $L$ by $L^{\vee}$. For an irreducible curve $D$ which is not necessarily smooth, we denote by $P_a(D)$ the arithmetic genus of $D$. For irreducible divisors $D_1$ and $D_2$ on a surface $X$, the arithmetic genus of $D_1+D_2$ is denoted as $P_a(D_1)+P_a(D_2)+D_1.D_2-1$. By induction, the arithmetic genus of a non-zero effective divisor $D$ on $X$ which is not irreducible is also defined. We denote it by the same notation $P_a(D)$. It follows that $2P_a(D)-2=D.(K_X+D)$ by the adjunction formula.

Let $C$ be a curve. Then the gonality of $C$ is the minimal degree of pencils on $C$. It is well known that if $C$ is a smooth plane curve of degree $d\geq5$, the gonality of $C$ is $d-1$.

Let $X$ be a surface. Then we denote the Picard lattice of $X$ by $\Pic(X)$, and call the rank of it the Picard number of $X$. If the Picard number of $X$ is $\rho$, then by the Hodge index theorem, the signature of $\Pic(X)$ is $(1,\rho-1)$. Note that this implies that $D_1^2D_2^2\leq(D_1.D_2)^2$ for two divisors $D_1$ and $D_2$ on $X$ satisfying  $D_1^2>0$ and $D_2^2>0$. 

Let $X$ be a smooth hypersurface of degree $d$ in $\mathbb{P}^3$. For a hyperplane section $C$ of $X$, we denote the class of it in $\Pic(X)$ by $H$. For an integer $l$, $H^{\otimes l}$ is often denoted as $\mathcal{O}_X(l)$. By the adjunction formula, $K_X\cong\mathcal{O}_X(d-4)$. For a vector bundle $E$ on $X$, we will write $E\otimes\mathcal{O}_X(l)=E(l)$.

\section{Preliminaries}

Let $X$ be a smooth quintic hypersurface in $\mathbb{P}^3$. In this section, we recall several fundamental notions concerning line bundles on $X$, and give some useful propositions about them.  Let $D$ be a divisor on $X$ and let $C$ be a smooth hyperplane section of $X$. First of all, since $K_X\cong\mathcal{O}_X(1)$, the Riemann-Roch theorem for $\mathcal{O}_X(D)$ is described as follows.
$$\chi(\mathcal{O}_X(D))=\displaystyle\frac{1}{2}D.(D-C)+\chi(\mathcal{O}_X),$$
where $\chi(\mathcal{O}_X(D))=h^0(\mathcal{O}_X(D))-h^1(\mathcal{O}_X(D))+h^2(\mathcal{O}_X(D))$. Note that since $h^0(K_X)=h^0(\mathcal{O}_{\mathbb{P}^3}(1))=4$, $h^1(\mathcal{O}_X)=0$, and $h^0(\mathcal{O}_X)=1$, we have $\chi(\mathcal{O}_X)=5$. 

\smallskip

\smallskip

\noindent The Serre duality for $\mathcal{O}_X(D)$ is given by 
$$h^i(\mathcal{O}_X(D))=h^{2-i}(\mathcal{O}_X(C-D))\;(0\leq i\leq 2).$$
By the Riemann-Roch theorem, if $h^0(\mathcal{O}_X(C-D))=0$, $h^0(\mathcal{O}_X(D))\geq\chi(\mathcal{O}_X(D))$. Hence, the following assertion is useful for estimating the value of $h^0(\mathcal{O}_X(D))$.

\newtheorem{prop}{Proposition}[section]

\begin{prop} Let $k$ be an integer with $0\leq k\leq 4$, and let $D$ be a non-zero effective divisor on $X$ such that $P_a(D)=C.D+1-k$. If $C.D\geq 7-k$, then $h^0(\mathcal{O}_X(C-D))=0$. \end{prop}

\noindent{\bf{Remark 2.1}}. The linear system $|\mathcal{O}_C(1)|$ gives an embedding $C\hookrightarrow\mathbb{P}^2$. Since $C$ is a plane quintic, $C$ is a 4-gonal curve. Therefore, if $L$ is a line bundle on $C$ satisfying $h^0(L)\geq2$, then $\deg(L)\geq h^0(L)+2$.

\smallskip

\smallskip

\noindent If $D$ is a non-zero effective divisor, then the arithmetic genus $P_a(D)$ is given as follows.
$$P_a(D)=\displaystyle\frac{1}{2}D.(D+C)+1.$$
If $D$ is reduced and irreducible, then $P_a(D)\geq0$. Before the proof of Proposition 2.1, we prepare the following lemmas.

\newtheorem{lem}{Lemma}[section]

\begin{lem} Let $D$ be a divisor on $X$ satisfying $C.D=1$. Then the following conditions are equivalent.

\smallskip

\smallskip

\noindent {\rm{(a)}} $h^0(\mathcal{O}_X(D))>0$.

\smallskip

\noindent {\rm{(b)}} $h^0(\mathcal{O}_X(D))=1$, $h^0(\mathcal{O}_X(C-D))=2$, and $h^1(\mathcal{O}_X(D))=0$.

\smallskip

\noindent {\rm{(c)}} $D^2=-3$.\end{lem}

{\it{Proof}}. (a) $\Longrightarrow$ (b). By the hypothesis, we may assume that $D$ is effective. Since $C.D=1$, by the ampleness of $C$, $D$ is reduced and irreducible. Therefore, $P_a(D)\geq0$, and hence, we have $D^2\geq-3$. This means that 
$$h^0(\mathcal{O}_X(D))+h^0(\mathcal{O}_X(C-D))\geq\chi(\mathcal{O}_X(D))\geq 3.\leqno (2.1)$$

On the other hand, since $C.D=1$ and $C.(C-D)=4$, by Remark 2.1, we have $h^0(\mathcal{O}_C(D))=1$ and $h^0(\mathcal{O}_C(C-D))\leq 2$. Since $C.(D-C)=-4$, by the ampleness of $C$, we have $h^0(\mathcal{O}_X(D-C))=0$. Since $h^0(\mathcal{O}_X(-D))=0$, by the exact sequence
$$0\longrightarrow \mathcal{O}_X(-D)\longrightarrow \mathcal{O}_X(C-D)\longrightarrow\mathcal{O}_C(C-D)\longrightarrow0,$$
and
$$0\longrightarrow \mathcal{O}_X(D-C)\longrightarrow \mathcal{O}_X(D)\longrightarrow\mathcal{O}_C(D)\longrightarrow0,$$
we have $h^0(\mathcal{O}_X(C-D))\leq2$ and $h^0(\mathcal{O}_X(D))=1$. By the inequality (2.1), we have $h^0(\mathcal{O}_X(C-D))=2$ and  $\chi(\mathcal{O}_X(D))=3$. Hence, we have $h^1(\mathcal{O}_X(D))=0$.

\smallskip

\smallskip

\noindent (b) $\Longrightarrow$ (c). Since $\chi(\mathcal{O}_X(D))=3$ and $C.D=1$, we have $D^2=-3$. 

\smallskip

\smallskip

\noindent (c) $\Longrightarrow$ (a). Since $C.D=1$, $h^0(\mathcal{O}_X(D))+h^0(\mathcal{O}_X(C-D))\geq\chi(\mathcal{O}_X(D))=3$. If $h^0(\mathcal{O}_X(D))=0$, we have $h^0(\mathcal{O}_X(C-D))\geq3$. Since $C.(C-D)=4$, by Remark 2.1, we have $h^0(\mathcal{O}_C(C-D))\leq2$. Since $-D.C=-1$, $h^0(\mathcal{O}_X(-D))=0$. Hence, by the exact sequence
$$0\longrightarrow \mathcal{O}_X(-D)\longrightarrow \mathcal{O}_X(C-D)\longrightarrow\mathcal{O}_C(C-D)\longrightarrow0,$$
we have $h^0(\mathcal{O}_X(C-D))\leq h^0(\mathcal{O}_C(C-D))$. This is a contradiction. Therefore, we have $h^0(\mathcal{O}_X(D))>0$.$\hfill\square$

\smallskip

\smallskip

\begin{lem} Let $D$ be an effective divisor on $X$ with $C.D=2$. $D^2\leq-6$ if and only if one of the following cases occurs.

\smallskip

\smallskip

\noindent {\rm(a)} There exists a curve $D_1$ on $X$ with $D=2D_1$, $D_1^2=-3$, and $C.D_1=1$.

\smallskip

\noindent {\rm(b)} There exist curves $D_1$ and $D_2$ with $D=D_1+D_2$, $D_1.D_2=0$, $D_i^2=-3$, and $C.D_i=1\;(i=1,2)$.\end{lem}

{\it{Proof}}. Assume that $D^2\leq-6$. If $D$ is reduced and irreducible, then $P_a(D)\geq0$. This means that $D^2\geq-4$. Hence, there exists a non-trivial effective decomposition $D=D_1+D_2$. Since $C.D=2$, we have $C.D_1=C.D_2=1$. Hence, $D_1$ and $D_2$ are reduced and irreducible and, by Lemma 2.1, $D_i^2=-3\;(i=1,2)$. If $D_1=D_2$, then $D=2D_1$ and  $D^2=-12$. If $D_1\neq D_2$, then $D_1.D_2\geq0$. Hence, we have $D^2\geq -6$. By the assumption, we have $D^2=-6$. Therefore, we have $D_1.D_2=0$. The converse assertion is clear. $\hfill\square$

$\;$

{\it{Proof of Proposition 2.1}}. If $C.D\geq\max\{6,7-k\}$, then $C.(C-D)\leq-1$. Hence, $h^0(\mathcal{O}_X(C-D))=0$. We consider the case where $2\leq k\leq 4$ and $C.D=5$. Assume that $|C-D|\neq\emptyset$. Then, by the ampleness of $C$, we have $\mathcal{O}_X(D)\cong\mathcal{O}_X(C)$. This contradicts the assumption that $P_a(D)=C.D+1-k$. Assume that $3\leq k\leq 4$ and $C.D=4$. Then $C.(C-D)=1$ and $(C-D)^2=1-2k\leq -5$. By Lemma 2.1, we have $|C-D|=\emptyset$. We consider the case where $k=4$ and $C.D=3$. Assume that $|C-D|\neq\emptyset$ and let $\Gamma$ be a member of $|C-D|$. Since $\Gamma.C=2$ and $\Gamma^2=-6$, by Lemma 2.2, there exist curves $D_1$ and $D_2$ satisfying the conditions that $\Gamma=D_1+D_2$, $D_1.D_2=0$, $D_i^2=-3$, and $C.D_i=1\;(i=1,2)$. Since $|C-\Gamma|=|D|\neq\emptyset$, there exists a hyperplane in $\mathbb{P}^3$ containing $\Gamma$. This is a contradiction. Hence, we have the assertion.$\hfill\square$

\smallskip

\smallskip

\begin{prop} Let $D$ be a non-zero effective divisor on $X$. If $h^1(\mathcal{O}_X(-D))=0$, then $P_a(D)\geq0$. \end{prop}

{\it{Proof}}. By the Serre duality, $h^1(\mathcal{O}_X(C+D))=0$, and $h^2(\mathcal{O}_X(C+D))=0$. Hence, by the Riemann-Roch theorem, we have $h^0(\mathcal{O}_X(C+D))=P_a(D)+4$. Since $h^0(\mathcal{O}_X(C+D))\geq h^0(\mathcal{O}_X(C))=4$, we have $P_a(D)\geq0$. $\hfill\square$

$\;$

\noindent By Proposition 2.2, for any effective divisor $D$ satisfying the condition (a) or (b) as in Lemma 2.2, we have $h^1(\mathcal{O}_X(-D))\neq0$. In general, the vanishing condition of the cohomology of the sheaf as in Proposition 2.2 can be characterized by the following notion.

\newtheorem{df}{Definition}[section]

\begin{df} Let $m$ be a positive integer. Then a non-zero effective divisor $D$ on $X$ is called {\rm{$m$-connected}} if $D_1.D_2\geq m$, for each effective decomposition $D=D_1+D_2$.\end{df}
\noindent If a non-zero effective divisor $D$ is 1-connected, then $h^0(\mathcal{O}_D)=1$. Therefore, by the exact sequence
$$0\longrightarrow \mathcal{O}_X(-D)\longrightarrow \mathcal{O}_X\longrightarrow\mathcal{O}_D\longrightarrow 0,$$
we have $h^1(\mathcal{O}_X(-D))=0$. Hence, by Proposition 2.1, we have $P_a(D)\geq0$. An effective divisor $D$ satisfying the condition as in Lemma 2.2 (b) is not 1-connected, but reduced. Hence, such a divisor $D$ is not contained by any hyperplane in $\mathbb{P}^3$. On the other hand, since it satisfies $h^0(\mathcal{O}_X(D-C))=0$, by the Riemann-Roch theorem and Remark 2.1, the conditions $h^0(\mathcal{O}_X(D))=1$ and $h^1(\mathcal{O}_X(D))=0$ are also satisfied. Conversely, any non-zero effective divisor $D$ which is not 1-connected is characterized as follows, under the condition that $\mathcal{O}_X(D)$ is initialized and $h^1(\mathcal{O}_X(D))=0$. 

\begin{prop} Let $D$ be a non-zero effective divisor. If $h^0(\mathcal{O}_X(D-C))=0$ and $h^1(\mathcal{O}_X(D))=0$, then $h^1(\mathcal{O}_X(-D))=0$ or there exist smooth rational curves $D_1$ and $D_2$ on $X$ such that $D=D_1+D_2$, $D_1.D_2=0$, and $C.D_i=1\;(i=1,2)$.\end{prop}

\noindent {\it{Proof}}. By the ampleness of $C$, $C.D\geq1$. If $C.D=1$, then $D$ is reduced and irreducible. Hence, $h^1(\mathcal{O}_X(-D))=0$. Assume that $C.D\geq 2$. If $C.D\geq 6$, then $C.(C-D)\leq-1$ and hence, $h^0(\mathcal{O}_C(C-D))=0$. Since $h^1(\mathcal{O}_X(C-D))=0$, by the exact sequence
$$0\longrightarrow\mathcal{O}_X(-D)\longrightarrow\mathcal{O}_X(C-D)\longrightarrow\mathcal{O}_C(C-D)\longrightarrow0,\leqno (2.2)$$
we have $h^1(\mathcal{O}_X(-D))=0$.  Assume that $C.D\leq5$. Since $C.D\geq2$, we have $0\leq C.(C-D)\leq 3$. By Remark 2.1, we have $h^0(\mathcal{O}_C(C-D))\leq1$. By the exact sequence (2.2), we have $h^0(\mathcal{O}_X(C-D))\leq1$. If $h^0(\mathcal{O}_X(C-D))=1$, then $h^0(\mathcal{O}_C(C-D))=1$ and  $h^1(\mathcal{O}_X(-D))=0$. Hence, it is sufficient to consider the case where $h^0(\mathcal{O}_X(C-D))=0$. By Remark 2.1, we have $h^0(\mathcal{O}_C(D))\leq 3$. By the assumption that $h^0(\mathcal{O}_X(D-C))=0$ and the exact sequence
$$0\longrightarrow\mathcal{O}_X(D-C)\longrightarrow\mathcal{O}_X(D)\longrightarrow\mathcal{O}_C(D)\longrightarrow0,$$
we have $h^0(\mathcal{O}_X(D))\leq 3$. If $h^0(\mathcal{O}_X(C-D))=0$, by the assumption that $h^1(\mathcal{O}_X(D))=0$, we have $h^0(\mathcal{O}_X(D))=\chi(\mathcal{O}_X(D))$, and hence, $-8\leq D^2-C.D\leq -4$. Therefore, it is sufficient to show the following lemma.

\begin{lem} Let $D$ be a non-zero effective divisor with $2\leq C.D\leq 5$ satisfying the condition that $h^1(\mathcal{O}_X(D))=0$, $h^0(\mathcal{O}_X(D-C))=0$, and $h^0(\mathcal{O}_X(C-D))=0$. Moreover, we assume that there exists an integer $k$ with $2\leq k\leq 4$ satisfying $D^2=C.D-2k$. If $D$ is not 1-connected, then there exist smooth rational curves $D_1$ and $D_2$ on $X$ such that $D=D_1+D_2$, $D_1.D_2=0$, and $C.D_i=1\;(i=1,2)$.\end{lem}

{\it{Proof}}. Assume that $D$ is not 1-connected, and let $D=D_1+D_2$ be a non-trivial effective decomposition with $D_1.D_2\leq0$. Since $D^2=D.D_1+D.D_2$, we may assume that $D.D_1=D^2-D.D_2\leq\displaystyle\frac{D^2}{2}$. Then we have
$$\displaystyle\frac{D^2}{2}\leq D.D_2=D_2^2+D_1.D_2\leq D_2^2.\leqno (2.3)$$
By the assumption that $D^2=C.D-2k$ and $h^0(\mathcal{O}_X(C-D))=0$ and the inequality (2.3), we have
$$\chi(\mathcal{O}_X(D_2))\geq h^0(\mathcal{O}_X(D))+\displaystyle\frac{1}{4}(C.D_1-C.D_2+2k).\leqno (2.4)$$

Assume that $C.D_1<C.D_2$. Since $C.D_1\geq1$, we have $C.D_2\geq2$. Hence, we have $C.(C-D_2)\leq 3$. By Remark 2.1, we have $h^0(\mathcal{O}_X(C-D_2))\leq1$. Assume that $h^0(\mathcal{O}_X(C-D_2))=0$. Since $C.D\leq 5$ and $2\leq k\leq 4$, we have
$$C.D_1-C.D_2+2k=2C.D_1-C.D+2k\geq 2k+2-C.D\geq 2k-3\geq 1.$$
Hence, by the inequality (2.4), we have $h^0(\mathcal{O}_X(D_2))>h^0(\mathcal{O}_X(D))$. This contradicts the fact that $\mathcal{O}_X(D_2)\subset\mathcal{O}_X(D)$. 

Assume that $h^0(\mathcal{O}_X(C-D_2))=1$. Since $h^0(\mathcal{O}_X(D-C))=0$, we have $C.(C-D_2)\geq1$. Since $C.D_1\geq1$ and $C.D_2\leq 4$, if $k=4$, then $C.D_1-C.D_2+2k-4\geq1$. By the inequality (2.4), we have $h^0(\mathcal{O}_X(D_2))>h^0(\mathcal{O}_X(D))$. However, by the same reason as above, this is a contradiction. If $k=2$ or 3, $C.D_1-C.D_2+2k-4\geq-3$. By the inequality (2.4), we have 
$$h^0(\mathcal{O}_X(D_2))\geq h^0(\mathcal{O}_X(D))=\chi(\mathcal{O}_X(D))=5-k\geq 2.$$
By the assumption that $h^0(\mathcal{O}_X(D-C))=0$, we have $h^0(\mathcal{O}_X(D_2-C))=0$. By Remark 2.1, we have $h^0(\mathcal{O}_C(D_2))\leq2$. Hence, by the exact sequence
$$0\longrightarrow\mathcal{O}_X(D_2-C)\longrightarrow\mathcal{O}_X(D_2)\longrightarrow\mathcal{O}_C(D_2)\longrightarrow0,\leqno(2.5)$$
we have $k=3$, $h^0(\mathcal{O}_C(D_2))=2$, and $C.D_2=4$. Since $C.(C-D_2)=1$, by Lemma 2.1, we have $(C-D_2)^2=-3$ and hence, $D_2^2=0$. Since $C.D\leq5$, we have $C.D_1=1$. By the same reason, $D_1^2=-3$. Since $D^2=C.D-6=-1$, we have $D_1.D_2=1$. This contradicts the hypothesis that $D_1.D_2\leq0$. By the above argument, we have $C.D_1\geq C.D_2$. 

Since $C.(C-D_2)\leq4$, by Remark 2.1, we have $h^0(\mathcal{O}_X(C-D_2))\leq 2$. If $h^0(\mathcal{O}_X(C-D_2))=0$, by the inequality (2.4), $h^0(\mathcal{O}_X(D_2))>h^0(\mathcal{O}_X(D))$. Hence, by the same reason as above, we have a contradiction.

Assume that $h^0(\mathcal{O}_X(C-D_2))=1$. Since $C.D_1\geq C.D_2$, if $3\leq k\leq 4$, then, by the inequality (2.4), we have $h^0(\mathcal{O}_X(D_2))>h^0(\mathcal{O}_X(D))$. By the same reason as above, we have a contradiction. Assume that $k=2$. Since $h^0(\mathcal{O}_X(D))=3$, by the inequality (2.4), we have $h^0(\mathcal{O}_X(D_2))\geq 3$. Since $|C-D_2|\neq\emptyset$ and $h^0(\mathcal{O}_X(D-C))=0$, we have $C.(C-D_2)\geq 1$. Since $C.D_2\leq 4$, by Remark 2.1 and the exact sequence (2.5), we have a contradiction.

Assume that $h^0(\mathcal{O}_X(C-D_2))=2$. Then, we have $C.(C-D_2)=4$. Hence, we have $C.D_2=1$. By Lemma 2.1, we have $D_2^2=-3$. Hence, by the inequality (2.3), we have $C.D\leq 2k-6$. Since $C.D\geq2$, by the assumption that $2\leq k\leq 4$, we have $k=4$, $C.D=2$, $D_1.D_2=0$ and $D^2=-6$. Since $C.D_1=C.D_2=1$, by Lemma 2.1, we have $P_a(D_1)=P_a(D_2)=0$. Therefore,  we have the assertion of Lemma 2.3. $\hfill\square$

\smallskip

\smallskip

\noindent We can construct an example of an effective divisor $D$ on a smooth quintic hypersurface $X$ satisfying the condition as in Lemma 2.2 (b).

$\;$

\noindent {\bf{Example 2.1}}. Let $X$ be the quintic hypersurface in $\mathbb{P}^3$ defined by the equation $x_0^5+x_1^5+x_2^5+x_3^5=0$ for a suitable homogeneous coordinate $(x_0:x_1:x_2:x_3)$ on $\mathbb{P}^3$. Let $D_1$ and $D_2$ be lines on $X$ which are defined by the equations $x_0+x_1=x_2+x_3=0$ and $x_0+x_2=x_1+\omega x_3=0$, respectively. Here, $\omega$ is a primitive 5 th root of 1. If we let $C$ be the hyperplane section of $X$ defined by the equation $x_1^5+x_2^5+x_3^5=0$, then $C.D_1=C.D_2=1$. Obviously, we have $D_1.D_2=0$. Hence, the divisor $D=D_1+D_2$ on $X$ is not 1-connected.

$\;$

\noindent By Proposition 2.3, we have the following assertion.

$\;$

\noindent{\bf{Corollary 2.1}}. Let $D$ be a non-zero effective divisor satisfying the condition that  $h^1(\mathcal{O}_X(D))=0$ and $h^0(\mathcal{O}_X(D-C))=0$. If $D^2>-6$, then we have $h^1(\mathcal{O}_X(-D))=0$.

\section{ACM bundles on quintic hypersurfaces in $\mathbb{P}^3$}

Let $X$ be as in section 2. In this section, we recall a well known fact about aCM bundles on $X$ and prepare a proposition to prove our main theorem. 

\begin{df} We call a vector bundle $E$ on $X$ an {\rm{arithmetically Cohen-Macaulay}} {\rm{(aCM}} for short{\rm{)}} bundle if $H^1(E(l))=0$ for all integers $l\in\mathbb{Z}$.
\end{df}

\begin{df} We say that a vector bundle $E$ on $X$ is {\rm{initialized}} if it satisfies the conditions $H^0(E)\neq0$ and $H^0(E(-1))=0$.

\end{df}

\noindent For an aCM bundle $E$ on $X$, we consider the graded module $H_{\ast}^0(E):=\bigoplus_{l\in\mathbb{Z}}H^0(E(l))$ over the homogeneous coordinate ring of $X$.  First of all, we mention the following result concerning the minimal number of generators of it.

\begin{prop}{\rm{([2, Theorem 3.1 and Corollary 3.5])}}. Let $E$ be an aCM bundle of rank $r$ on $X$, and let $\mu(E)$ be the minimal number of generators of $H_{\ast}^0(E)$. Then we get $\mu(E)\leq 5r$. Moreover, if $E$ is initialized, then $h^0(E)\leq5r$.\end{prop}

\noindent An initialized aCM bundle $E$ of rank $r$ with $h^0(E)=5r$ is called an {\it{Ulrich}} bundle. An Ulrich bundle of rank $r$ is characterized as an initialized aCM bundle whose Hilbert polynomial is equal to $5r\binom{t+2}{2}$. In particular, the line bundle $\mathcal{O}_X(D)$ defined by an effective divisor $D$ on $X$ with $D^2=D.K_X=10$ as in Theorem 1.1 (ii) is Ulrich. If there exists such a line bundle $\mathcal{O}_X(D)$ on $X$, by taking the minimal free resolution of $H_{\ast}^0(\mathcal{O}_X(D))$ as a module over the homogeneous coordinate ring of $\mathbb{P}^3$, it follows that $X$ is linear determinantal (i.e., $X$ is defined as the zero locus of the determinant of a $5\times5$-matrix of linear forms). 

$\;$

\noindent{\bf{Remark 3.1}}. A line bundle on $X$ is aCM if and only if the dual of it is aCM.

\smallskip

\smallskip

\begin{prop} Let $D$ be a non-zero effective divisor on $X$, let $C$ be a smooth hyeperplane section of $X$, and let $k$ be a positive integer satisfying $C.D+5<5k$. If $h^1(\mathcal{O}_X(lC-D))=0$ for $0\leq l\leq k$, then $\mathcal{O}_X(D)$ is aCM. \end{prop}

\smallskip

\smallskip

\noindent {\bf{Proof}}. Since $K_C=\mathcal{O}_C(2)$, by the Serre duality, for $n\geq1$, we have 
$$h^1(\mathcal{O}_C((n+1)C-D))=h^0(\mathcal{O}_C(D-(n-1)C)).$$
If $n\geq k$, by the assumption, $\deg(\mathcal{O}_C(D-(n-1)C))=C.D-5(n-1)<0$, and hence, we have $h^1(\mathcal{O}_C((n+1)C-D))=0$. Since $h^1(\mathcal{O}_X(lC-D))=0$ for $1\leq l\leq k$, by the exact sequence
$$0\longrightarrow\mathcal{O}_X(nC-D)\longrightarrow\mathcal{O}_X((n+1)C-D)\longrightarrow\mathcal{O}_C((n+1)C-D)\longrightarrow0,$$
for any integer $n\geq1$, we have $h^1(\mathcal{O}_X(nC-D))=0$ by induction.

On the other hand, if $m\geq0$, we have $h^0(\mathcal{O}_C(-mC-D))=0$. Since $h^1(\mathcal{O}_X(-D))=0$, by the exact sequence
$$0\longrightarrow\mathcal{O}_X(-(m+1)C-D)\longrightarrow\mathcal{O}_X(-mC-D)\longrightarrow\mathcal{O}_C(-mC-D)\longrightarrow0,$$
for any integer $m\geq0$, we have $h^1(\mathcal{O}_X(-mC-D))=0$ by induction. Hence, $\mathcal{O}_X(-D)$ is aCM. By Remark 3.1, $\mathcal{O}_X(D)$ is also aCM. $\hfill\square$

\section{Proof of Theorem 1.1}

Let $X$ be a smooth quintic hypersurface in $\mathbb{P}^3$, and let $C$ be a smooth hyperplane section of $X$. Let $D$ be a non-zero effective divisor of arithmetic genus $P_a(D)$. In this section, we give a proof of our main theorem. First of all, we have the following assertion.

\begin{prop} If $\mathcal{O}_X(D)$ is aCM and initialized, then the following conditions are satisfied.

\smallskip

\smallskip

\noindent {\rm{(a)}} $C.D-3\leq P_a(D)\leq C.D+1$.

\smallskip

\noindent {\rm{(b)}} If $P_a(D)=C.D-3$, then $3\leq C.D\leq 6$.

\smallskip

\noindent {\rm{(c)}} If $P_a(D)=C.D-2$, then $2\leq C.D\leq 7$.

\smallskip

\noindent {\rm{(d)}} If $P_a(D)=C.D-1$, then $C.D=1$ or $4\leq C.D\leq 8$.

\smallskip

\noindent {\rm{(e)}} If $P_a(D)=C.D$, then $6\leq C.D\leq 9$.

\smallskip

\noindent {\rm{(f)}} If $P_a(D)=C.D+1$, then $7\leq C.D\leq 10$. \end{prop}

\noindent By the ampleness of $C$, we have $C.(C-D)\leq 4$. Hence, $h^0(\mathcal{O}_C(C-D))\leq2$, by Remark 2.1. By the exact sequence
$$0\longrightarrow\mathcal{O}_X(-D)\longrightarrow\mathcal{O}_X(C-D)\longrightarrow\mathcal{O}_C(C-D)\longrightarrow0,\leqno (4.1)$$
we have  $h^0(\mathcal{O}_X(C-D))\leq2$. Hence, we divide the assertion of Proposition 4.1 into the following three lemmas depending on the value of $h^0(\mathcal{O}_X(C-D))$.

\begin{lem} Assume that $h^0(\mathcal{O}_X(C-D))=2$. Then $P_a(D)=0$ and $C.D=1$. \end{lem}

{\it{Proof}}. By the exact sequence (4.1), we have $h^0(\mathcal{O}_C(C-D))=2$. Hence, we have $C.D=1$. By Lemma 2.1, we have $D^2=-3$ and hence, $P_a(D)=0$. $\hfill\square$

\begin{lem} Assume that $h^0(\mathcal{O}_X(C-D))=1$. If $\mathcal{O}_X(D)$ is aCM and initialized, then one of the following cases occurs.

\smallskip

\smallskip

\noindent {\rm{(i)}} $P_a(D)=0$ and $C.D=2$.

\smallskip

\noindent {\rm{(ii)}} $P_a(D)=1$ and $C.D=3$.

\smallskip

\noindent {\rm{(iii)}} $P_a(D)=3$ and $C.D=4$. \end{lem}

{\it{Proof}}. Since $h^1(\mathcal{O}_X(D))=0$, by the Riemann-Roch theorem, we have 
$$h^0(\mathcal{O}_X(D))=P_a(D)+3-C.D.\leqno (4.2)$$
Since $|C-D|\neq\emptyset$, $C.(C-D)\geq0$. Hence, we have $C.D\leq 5$. If $C.D=5$, $\mathcal{O}_X(D)\cong\mathcal{O}_X(C)$. This contradicts the assumption that $\mathcal{O}_X(D)$ is initialized. Hence, we have $C.D\leq4$. Since $\mathcal{O}_X(D)$ is aCM and initialized, by the exact sequence
$$0\longrightarrow\mathcal{O}_X(D-C)\longrightarrow\mathcal{O}_X(D)\longrightarrow\mathcal{O}_C(D)\longrightarrow0,$$
we have $h^0(\mathcal{O}_C(D))=h^0(\mathcal{O}_X(D)).$ By Remark 2.1, $h^0(\mathcal{O}_X(D))\leq2$. Assume that $h^0(\mathcal{O}_X(D))=2$. Since we have $C.D=4$, by the equality (4.2), $P_a(D)=3$. Assume that $h^0(\mathcal{O}_X(D))=1$. Then $P_a(D)=C.D-2\leq 2$. If $P_a(D)=2$, then $C.D=4$. Hence, $C.(C-D)=1$ and $(C-D)^2=-5$. However, since $|C-D|\neq\emptyset$, this contradicts Lemma 2.1. Since $h^1(\mathcal{O}_X(-D))=0$, by Proposition 2.2, $P_a(D)\geq0$. Hence, we have the assertion. $\hfill\square$

\begin{lem} Assume that $h^0(\mathcal{O}_X(C-D))=0$. If $\mathcal{O}_X(D)$ is aCM and initialized, then one of the following cases occurs.

\smallskip

\smallskip

\noindent {\rm{(i)}} $P_a(D)=C.D-3$ and $3\leq C.D\leq 6$.

\smallskip

\noindent {\rm{(ii)}} $P_a(D)=C.D-2$ and $4\leq C.D\leq 7$.

\smallskip

\noindent {\rm{(iii)}} $P_a(D)=C.D-1$ and $5\leq C.D\leq 8$.

\smallskip

\noindent {\rm{(iv)}} $P_a(D)=C.D$ and $6\leq C.D\leq 9$.

\smallskip

\noindent {\rm{(v)}} $P_a(D)=C.D+1$ and $7\leq C.D\leq 10$. \end{lem}

{\it{Proof}}. First of all, by Proposition 3.1, we have $h^0(\mathcal{O}_X(D))\leq 5$. Let $k$ be an integer with $0\leq k\leq 4$ such that $h^0(\mathcal{O}_X(D))=5-k$. Since $h^1(\mathcal{O}_X(D))=0$, by the Riemann-Roch theorem, we have 
$$P_a(D)=C.D+1-k.\leqno (4.3)$$

On the other hand, we have $h^0(\mathcal{O}_X(D))+5-C.D=\chi(\mathcal{O}_X(2C-D))$. Since $h^1(\mathcal{O}_X(2C-D))=0$, we have 
$$C.D\leq h^0(\mathcal{O}_X(D))+5.\leqno (4.4)$$

If $h^0(\mathcal{O}_X(D))=1$, then we have $P_a(D)=C.D-3$ by the equality (4.3). Since $h^1(\mathcal{O}_X(-D))=0$, by Proposition 2.2, we have $P_a(D)\geq0$. Hence, we have $C.D\geq3$. By the inequality (4.4), we have $C.D\leq6$. 

Assume that $h^0(\mathcal{O}_X(D))\geq2$. Then $0\leq k\leq 3$. By the same argument as in the proof of Lemma 4.2, we have $h^0(\mathcal{O}_C(D))=5-k$. By Remark 2.1, we have $C.D\geq 7-k$. By the inequality (4.4), we have $C.D\leq 10-k$. Hence, we have the assertion.$\hfill\square$

\smallskip

\smallskip

\noindent By Proposition 4.1, we have the following necessary condition for $\mathcal{O}_X(D)$ to be aCM and initialized.

\begin{prop} Let $k=C.D+1-P_a(D)$. If $\mathcal{O}_X(D)$ is aCM and initialized, then the following conditions are satisfied.

$\;$

\noindent {\rm{(i)}} $0\leq k\leq 4$.

\smallskip

\smallskip

\noindent {\rm{(ii)}} If $0\leq k\leq 1$, then $C.D=10-k$ and $h^0(\mathcal{O}_C(D-C))=0$. 

\smallskip

\smallskip

\noindent {\rm{(iii)}} If $k=2$, then the following conditions are satisfied.

\smallskip

\smallskip

{\rm{(a)}} $C.D=1$ or $4\leq C.D\leq 8$.

{\rm{(b)}} If $C.D=7$, then $h^0(\mathcal{O}_X(2C-D))=1$.

{\rm{(c)}} If $C.D=8$, then $h^0(\mathcal{O}_C(D-C))=0$ and $h^0(\mathcal{O}_C(D))=3$.

\smallskip

\smallskip

\noindent {\rm{(iv)}} If $3\leq k\leq 4$, then the following conditions are satisfied.

\smallskip

{\rm{(a)}} $k-1\leq C.D\leq 10-k$.

{\rm{(b)}} If $8-k\leq C.D\leq 10-k$, then $h^0(\mathcal{O}_C(D))=5-k$. 

\end{prop}

{\it{Proof}}. First of all, by Proposition 4.1, the assertion of (i) is clear. Assume that $C.D\geq 8-k$. By Proposition 2.1, we have $|C-D|=\emptyset$. Since $h^1(\mathcal{O}_X(D))=0$, by the Riemann-Roch theorem, we have $h^0(\mathcal{O}_X(D))=5-k$. Since $\mathcal{O}_X(D)$ is aCM and initialized, by the exact sequences
$$0\longrightarrow\mathcal{O}_X(D-C)\longrightarrow\mathcal{O}_X(D)\longrightarrow\mathcal{O}_C(D)\longrightarrow0,$$
and
$$0\longrightarrow\mathcal{O}_X(D-2C)\longrightarrow\mathcal{O}_X(D-C)\longrightarrow\mathcal{O}_C(D-C)\longrightarrow0,$$
we have $h^0(\mathcal{O}_C(D))=5-k$ and $h^0(\mathcal{O}_C(D-C))=0$. Therefore, by Proposition 4.1, we get the assertion of (iv). 

On the other hand, we have
$$h^0(\mathcal{O}_X(2C-D))=\chi(\mathcal{O}_X(2C-D))=10-k-C.D.\leqno (4.5)$$

(ii) We consider the case where $0\leq k\leq1$. By the above argument, it is sufficient to show that $C.D=10-k$. By Proposition 4.1, $7-k\leq C.D\leq 10-k$. Assume that $7-k\leq C.D\leq 8-k$. Then, by the equality (4.5), we have $h^0(\mathcal{O}_X(2C-D))\geq2$. On the other hand, since $|C-D|=\emptyset$ and  $h^1(\mathcal{O}_X(C-D))=0$, by the equality (4.5) and the exact sequence
$$0\longrightarrow\mathcal{O}_X(C-D)\longrightarrow\mathcal{O}_X(2C-D)\longrightarrow\mathcal{O}_C(2C-D)\longrightarrow0,$$
we have $h^0(\mathcal{O}_C(2C-D))=10-k-C.D$. By Remark 2.1, we have $C.(2C-D)\geq 12-k-C.D$. This contradicts the assumption that $0\leq k\leq1$. 

Assume that $C.D=9-k$. By the equality (4.5), $h^0(\mathcal{O}_X(2C-D))=1$. If $k=0$, we have $C.(2C-D)=1$ and $(2C-D)^2=-7$. However, this contradicts Lemma 2.1. Assume that $k=1$ and let $\Gamma$ be the member of $|2C-D|$. Then we have $P_a(\Gamma)=-1$. By Proposition 2.2, this means that $h^1(\mathcal{O}_X(D-2C))\neq0$. This contradicts the assumption that $\mathcal{O}_X(D)$ is aCM. Hence, we have the assertion.

\smallskip

\smallskip

(iii) If $k=2$ and $C.D=7$, then, by the equality (4.5), we have the assertion of (b). By Proposition 4.1 and the above argument, we have the assertion.  $\hfill\square$

$\;$

\noindent From now on we show that each condition from (ii) to (iv) as in Proposition 4.2 is a sufficient condition for $\mathcal{O}_X(D)$ to be aCM and initialized. Since the proof is long and complex, we divide the converse assertion of Proposition 4.2 into several propositions.

\begin{prop} Assume that $P_a(D)=C.D-3$. If the following conditions are satisfied, then $\mathcal{O}_X(D)$ is aCM and initialized. 

\smallskip

\smallskip

\noindent {\rm{(a)}} $3\leq C.D\leq 6$.

\smallskip

\noindent {\rm{(b)}} If $4\leq C.D\leq 6$, then $h^0(\mathcal{O}_C(D))=1$. \end{prop}

{\it{Proof}}. First of all, $\mathcal{O}_X(D)$ is initialized. Indeed, if $C.D\leq4$, then we have $C.(D-C)\leq-1$, and hence, $|D-C|=\emptyset$. Since $h^0(\mathcal{O}_C(1))=3$, if $5\leq C.D\leq 6$, by the assumption (b), $h^0(\mathcal{O}_C(D-C))=0$. Moreover, since $C.(D-2C)<0$, we have $h^0(\mathcal{O}_X(D-2C))=0$. By the exact sequence
$$0\longrightarrow\mathcal{O}_X(D-2C)\longrightarrow\mathcal{O}_X(D-C)\longrightarrow\mathcal{O}_C(D-C)\longrightarrow0,\leqno (4.6)$$
we have $h^0(\mathcal{O}_X(D-C))=0$. 

If $n\geq 3$, we have $C.D+5<5n$. Hence, by Proposition 3.2, it is sufficient to show that $h^1(\mathcal{O}_X(lC-D))=0$ for $0\leq l\leq 3$. By Proposition 2.1, we have $|C-D|=\emptyset$. If $C.D=3$, by Remark 2.1, we have $h^0(\mathcal{O}_C(D))=1$. Hence, by the assumption (b) and the exact sequence
$$0\longrightarrow\mathcal{O}_X(D-C)\longrightarrow\mathcal{O}_X(D)\longrightarrow\mathcal{O}_C(D)\longrightarrow0,\leqno (4.7)$$
we have $h^0(\mathcal{O}_X(D))=1$. Since $\chi(\mathcal{O}_X(D))=1$, we have 
$$h^1(\mathcal{O}_X(C-D))=h^1(\mathcal{O}_X(D))=0.$$
By the exact sequence (4.6) and (4.7), we have $h^1(\mathcal{O}_X(lC-D))=0$ for $2\leq l\leq3$.
Since $\mathcal{O}_X(D)$ is initialized, and satisfies $h^1(\mathcal{O}_X(D))=0$ and $P_a(D)\geq0$, by Corollary 2.1, we have $h^1(\mathcal{O}_X(-D))=0$. Hence, we have the assertion. $\hfill\square$

\begin{prop} Assume that $P_a(D)=C.D-2$. If the following conditions are satisfied, then $\mathcal{O}_X(D)$ is aCM and initialized. 

\smallskip

\smallskip

\noindent {\rm{(a)}} $2\leq C.D\leq 7$.

\smallskip

\noindent {\rm{(b)}} If $5\leq C.D\leq 7$, then $h^0(\mathcal{O}_C(D))=2$. \end{prop}

{\it{Proof}}. First of all, $\mathcal{O}_X(D)$ is initialized. Indeed, if $C.D\leq 4$, $C.(D-C)\leq-1$ and hence, we have $|D-C|=\emptyset$. Since $h^0(\mathcal{O}_C(1))=3$, if $5\leq C.D\leq7$, then the assumption (b) implies $h^0(\mathcal{O}_C(D-C))=0$. Moreover, since $C.(D-2C)<0$, we have $h^0(\mathcal{O}_X(D-2C))=0$. Hence, by the exact sequence (4.6), we have $h^0(\mathcal{O}_X(D-C))=0$. If $n\geq3$, then we have $C.D+5<5n$. Hence, by Proposition 3.2, it is sufficient to show that $h^1(\mathcal{O}_X(lC-D))=0$ for $0\leq l\leq3$. 

Assume that $2\leq C.D\leq 3$. By Remark 2.1 and the exact sequence (4.7), we have $h^0(\mathcal{O}_X(D))=h^0(\mathcal{O}_C(D))=1$. Since $\chi (\mathcal{O}_X(D))=2$, $h^0(\mathcal{O}_X(C-D))\geq1$. Since $2\leq C.(C-D)\leq 3$, we have $h^0(\mathcal{O}_C(C-D))\leq1$. Hence, by the exact sequence
$$0\longrightarrow\mathcal{O}_X(-D)\longrightarrow\mathcal{O}_X(C-D)\longrightarrow\mathcal{O}_C(C-D)\longrightarrow0,\leqno (4.8)$$
we have $h^0(\mathcal{O}_X(C-D))=h^0(\mathcal{O}_C(C-D))=1$. Hence, we have 
$$h^1(\mathcal{O}_X(C-D))=h^1(\mathcal{O}_X(D))=0.$$

Assume that $4\leq C.D\leq 7$. Since, by Proposition 2.1, $|C-D|=\emptyset$, we have $h^0(\mathcal{O}_X(D))\geq\chi(\mathcal{O}_X(D))=2$. If $C.D=4$, by Remark 2.1, $h^0(\mathcal{O}_C(D))\leq 2$. By the assumption (b) and the exact sequence (4.7), we have $h^0(\mathcal{O}_X(D))=h^0(\mathcal{O}_C(D))=2$. Hence, we have $h^1(\mathcal{O}_X(C-D))=h^1(\mathcal{O}_X(D))=0$. Hence, by the exact sequence (4.7), we have $h^1(\mathcal{O}_X(2C-D))=h^1(\mathcal{O}_X(D-C))=0$. Since $h^0(\mathcal{O}_C(D-C))=0$, by the exact sequence (4.6), we have 
$$h^1(\mathcal{O}_X(3C-D))=h^1(\mathcal{O}_X(D-2C))=0.$$
Since $\mathcal{O}_X(D)$ is initialized and satisfies $h^1(\mathcal{O}_X(D))=0$ and $P_a(D)\geq0$, by Corollary 2.1, we have $h^1(\mathcal{O}_X(-D))=0$. Therefore, the assertion follows. $\hfill\square$

\begin{prop} Assume that $P_a(D)=C.D-1$. If the following conditions are satisfied, then $\mathcal{O}_X(D)$ is aCM and initialized. 

\smallskip

\smallskip

\noindent {\rm{(a)}} $C.D=1$ or $4\leq C.D\leq 8$.

\noindent {\rm{(b)}} If $C.D=7$, then $h^0(\mathcal{O}_X(2C-D))=1$.

\noindent {\rm{(c)}} If $C.D=8$, then $h^0(\mathcal{O}_C(D-C))=0$ and $h^0(\mathcal{O}_C(D))=3$. \end{prop}

{\it{Proof}}. First of all, we show that $\mathcal{O}_X(D)$ is initialized. If $C.D=1$ or 4, then $C.(D-C)<0$, and hence, we have $|D-C|=\emptyset$. We consider the case where $5\leq C.D\leq7$. Assume that $|D-C|\neq\emptyset$. If $C.D=5$, then $\mathcal{O}_X(D)\cong\mathcal{O}_X(C)$. This contradicts the assumption  that $P_a(D)=4$. If $C.D=6$, then $(D-C)^2=-5$ and $C.(D-C)=1$. This contradicts Lemma 2.1. If $C.D=7$, then by the assumption (b), $|2C-D|\neq\emptyset$. Hence, if we let  $\Gamma$ be the member of $|D-C|$, $\Gamma$ is contained by a hyperplane in $\mathbb{P}^3$. Since $\Gamma.C=2$ and $\Gamma^2=-6$, by Lemma 2.2, this is a contradiction. Hence, $|D-C|=\emptyset$. If $C.D=8$, then $C.(D-2C)=-2$, and hence, $|D-2C|=\emptyset$. By the assumption (c) and the exact sequence (4.6), we have $h^0(\mathcal{O}_X(D-C))=0$. If $n\geq 3$, then $C.D+5<5n$. Hence, it is sufficient to show that $h^1(\mathcal{O}_X(lC-D))=0$ for $0\leq l\leq 3$. 

We show that $h^1(\mathcal{O}_X(C-D))=0$ and $h^0(\mathcal{O}_X(D))=h^0(\mathcal{O}_C(D))$. If $C.D=1$, by the exact sequence (4.7), we have $h^0(\mathcal{O}_X(D))=h^0(\mathcal{O}_C(D))=1$. By Lemma 2.1, we have $h^1(\mathcal{O}_X(C-D))=h^1(\mathcal{O}_X(D))=0$. 

Assume that $C.D=4$. Since $C.(C-D)=1$, we have $h^0(\mathcal{O}_C(C-D))\leq1$. By the exact sequence (4.8), we have $h^0(\mathcal{O}_X(C-D))\leq1$. Since $\chi(\mathcal{O}_X(D))=3$, we have $h^0(\mathcal{O}_X(D))\geq2$. By Remark 2.1 and the exact sequence (4.7), we have $h^0(\mathcal{O}_X(D))=h^0(\mathcal{O}_C(D))=2$. Hence, we have $h^0(\mathcal{O}_X(C-D))=1$ and $h^1(\mathcal{O}_X(C-D))=h^1(\mathcal{O}_X(D))=0$. 

We consider the case where $5\leq C.D\leq8$. By Proposition 2.1, we have $|C-D|=\emptyset$. Hence, $h^0(\mathcal{O}_X(D))\geq3$. We have $h^0(\mathcal{O}_X(D))=3$. Indeed, if $C.D=5$, then $h^0(\mathcal{O}_X(D))=h^0(\mathcal{O}_C(D))=3$, by Remark 2.1 and the exact sequence (4.7). Assume that $6\leq C.D\leq7$. Since $\mathcal{O}_X(D)$ is initialized, $h^0(\mathcal{O}_X(2C-D))\geq\chi(\mathcal{O}_X(2C-D))=8-C.D$. By the exact sequence
$$0\longrightarrow\mathcal{O}_X(C-D)\longrightarrow\mathcal{O}_X(2C-D)\longrightarrow\mathcal{O}_C(2C-D)\longrightarrow0,$$
we have $h^0(\mathcal{O}_C(2C-D))\geq8-C.D$. Since $C.(2C-D)=10-C.D$, by Remark 2.1, we have $h^0(\mathcal{O}_C(2C-D))=8-C.D$. Since $K_C\cong\mathcal{O}_C(2)$, by the Riemann-Roch theorem, we have $h^0(\mathcal{O}_C(D))=3$. By the exact sequence (4.7), we have $h^0(\mathcal{O}_X(D))=3$. If $C.D=8$, by the assumption (c) and the exact sequence (4.7), we have $h^0(\mathcal{O}_X(D))=3$. Hence, we have $h^1(\mathcal{O}_X(C-D))=h^1(\mathcal{O}_X(D))=0$. For each case as above, by using the exact sequence (4.7), we have $h^1(\mathcal{O}_X(2C-D))=h^1(\mathcal{O}_X(D-C))=0$. Since $P_a(D)\geq0$, by Corollary 2.1, we have $h^1(\mathcal{O}_X(-D))=0$. 

We show that $h^1(\mathcal{O}_X(3C-D))=0$. If $C.D=$1,4, or 8, then we have $h^0(\mathcal{O}_C(D-C))=0$. Hence, by the exact sequence (4.6), we have 
$$h^1(\mathcal{O}_X(3C-D))=h^1(\mathcal{O}_X(D-2C))=0.$$ 
Assume that $5\leq C.D\leq7$. Then $\chi(\mathcal{O}_X(2C-D))=8-C.D>0$. Since $\mathcal{O}_X(D)$ is initialized and $h^1(\mathcal{O}_X(2C-D))=0$, we have $h^0(\mathcal{O}_X(2C-D))>0$. Since $|C-D|=\emptyset$ and $(2C-D)^2\geq-5$, if we apply Corollary 2.1 to $\mathcal{O}_X(2C-D)$, we have $h^1(\mathcal{O}_X(3C-D))=h^1(\mathcal{O}_X(D-2C))=0$. Hence, the assertion follows. $\hfill\square$

\begin{prop} Let $k=C.D+1-P_a(D)$. Assume that $0\leq k\leq 1$ and $C.D=10-k$. If $h^0(\mathcal{O}_C(D-C))=0$, then $\mathcal{O}_X(D)$ is aCM and initialized. \end{prop}

{\it{Proof}}. First of all, we show that $\mathcal{O}_X(D)$ is initialized. First of all, we have $|D-2C|=\emptyset$. In fact, since $C.(D-2C)=-k$, if $|D-2C|\neq\emptyset$, then we have $k=0$ and $\mathcal{O}_X(D)\cong\mathcal{O}_X(2C)$. This contradicts the assumption that $P_a(D)=11$. By the exact sequence (4.6) and the assumption that $h^0(\mathcal{O}_C(D-C))=0$, we have $h^0(\mathcal{O}_X(D-C))=0$. If $n\geq4$, then we have $C.D+5<5n$. Hence, it is sufficient to show that $h^1(\mathcal{O}_X(lC-D))=0$ for $0\leq l\leq 4$. 

We show that $h^1(\mathcal{O}_X(C-D))=0$. Since $h^0(\mathcal{O}_C(D-C))=0$ and $0\leq k\leq 1$, we have $h^0(\mathcal{O}_C(2C-D))=0$. By the Riemann-Roch theorem, we have $h^0(\mathcal{O}_C(D))=5-k$. On the other hand, since $C.(C-D)=k-5<0$, we have $|C-D|=\emptyset$. Hence, we have $h^0(\mathcal{O}_X(D))\geq\chi(\mathcal{O}_X(D))=5-k$. By the exact sequence (4.7), we have $h^0(\mathcal{O}_X(D))=5-k$. Therefore, we have $h^1(\mathcal{O}_X(C-D))=h^1(\mathcal{O}_X(D))=0$. By the exact sequence (4.8), we have $h^1(\mathcal{O}_X(-D))=0$. By the exact sequence (4.6) and (4.7), we have $h^1(\mathcal{O}_X(lC-D))=0$ for $2\leq l\leq 3$. If $k=1$, then since $C.D+5<15$, by Proposition 3.2, the assertion is clear. Assume that $k=0$. Since $h^0(\mathcal{O}_C(D-C))=0$, we have $h^0(\mathcal{O}_C(D-2C))=0$ and hence, by the exact sequence
$$0\longrightarrow\mathcal{O}_X(D-3C)\longrightarrow\mathcal{O}_X(D-2C)\longrightarrow\mathcal{O}_C(D-2C)\longrightarrow0,$$
we have $h^1(\mathcal{O}_X(4C-D))=h^1(\mathcal{O}_X(D-3C))=0$. Hence, we have the assertion.$\hfill\square$

$\;$

\noindent {\bf{Acknowledgements}}

\smallskip

\smallskip

\noindent The author would like to thank the referee for some helpful comments.

\end{document}